\documentclass[12pt]{amsart}
\usepackage{amsmath,amsthm,latexsym,amscd,amsbsy,amssymb,url}
\usepackage[all]{xypic}
 \relax

\chardef\bslash=`\\ 

\makeatletter
\def\verbatim{\interlinepenalty\@M \@verbatim
  \leftskip\@totalleftmargin\advance\leftskip2pc
  \frenchspacing\@vobeyspaces \@xverbatim}
\makeatother
\newtheorem{thm}{Theorem}[section]
\newtheorem{cor}[thm]{Corollary}
\newtheorem{lem}[thm]{Lemma}

\newtheorem{ex}[thm]{Example}

\newtheorem{problem}[thm]{Problem}

\begin{document}


\title
{On Topological Groups of Monotonic Autohomeomorphisms}
\author{Raushan  Buzyakova}
\email{Raushan\_Buzyakova@yahoo.com}

\keywords{linearly ordered topological space, generalized ordered topological space, topology of point-wise convergence, paratopological group, topological group, monotonic autohomeomorphism}
\subjclass{54F05, 54H11, 54C35}


\begin{abstract}{
We study topological groups of monotonic autohomeomorphisms on a generalized ordered space $L$. We find a condition that is necessary and sufficient for the set of all monotonic autohomeomorphisms on $L$ along with the function composition and the topology of point-wise convergence to be  a topological group.
}
\end{abstract}

\maketitle
\markboth{R. Buzyakova}{On Topological Groups of Monotonic Autohomeomorphisms}
{ }

\section{Introduction}\label{S:introduction}
\par\bigskip\noindent
In this paper we study sets of monotonic autohomeomorphisms on generalized  ordered spaces endowed with the topology of point-wise convergence. Recall that a linearly ordered topological space, abbreviated as LOTS, is  a linearly ordered set along with the topology generated by sets in form $(a,b), [\min L, a), (a, \max L]$. A generalized ordered space, abbreviated as GO,  is a subspace of a linearly ordered space.  It is a result of \v Cech  that the topology of a generalized ordered space is generated by a collection of convex subsets \cite{BL}.  When we discuss several linearly ordered spaces, to distinguish their intervals we will use subscription as in $(a,b)_L$, where $a,b$ can be in $L$ or in  a  superspace understood from the context.
It is a traditional exercise that the set of all monotonic autohomeomorphisms $M(L)$ on a GO-space $L$ along with the operation of composition is a group. We observe that  this group along with the topology of point-wise convergence, denoted by $M_p(L)$,  is a paratopological group. Recall that a group $G$  along with a  topology on $G$ is a paratopological group if the group operation of $G$ is continuous with respect to the topology of $G\times G$. The operation of inversion, however, need not be continuous in $M_p(L)$. We, therefore, identify a condition that is necessary and sufficient for  $M_p(L)$ to be a topological group. In the main result of this work (Theorem \ref{thm:criterion}), we prove that given a GO-space $L$, the space $M_p(L)$ is a topological group if and only any set in from  $U(x; \{y\})=\{f: f(x) = y\}$ is open whenever $x$ is limit from at most one side. It is not hard to see that for a LOTS $L$, such sets are always open. Therefore, a corollary to our main result is  the recent result of B. Sorin \cite{Sor} that the group of order-preserving bijectiions of a linearly ordered space $L$  with the operation of composition and  endowed with the topology of point-wise convergence is a topological group. Sorin's argument uses the fact that the topology of point-wise convergence of the group of continuous extensions over the smallest linearly  ordered compactification is generated by sets dependent on points of $L$ only. This, however, is no longer true for a GO-space.  It is also worth mentioning  that the group of isometries on a metric space with the topology of point-wise convergence is a topological group too (see \cite[Theorem 3.5.1]{AT}). Since monotonic maps are either order-preserving or order-reversing, we can view them as the order counterparts of isometries on metric spaces. In general, the topology of point-wise convergence need not turn a group of autohomeomorphisms on a space (even linearly ordered space)  with operation of composition  into a topological group. It is easy to see that neither taking the inverse nor the operation of composition are continuous with respect to this topology even for autohomeomorphisms of $\mathbb Q$.  

\par\smallskip\noindent
Given a GO-space $L$, standard open sets of $M_p(L)$ are in form 
$$U=U(x_1,..,x_n; I_1,...,I_n)=\{f\in M(L): f(x_i)\in I_i, i=1,...,n\},$$
where $x_i$'s  are some fixed elements of $L$ and  $I_k$'s are open convex sets of $L$. Note that  $U(x_1,..,x_n; I_1,...,I_n)=U(x_1; I_1)\cap ...\cap U(x_n;I_n)$. An unordered pair of elements is a gap in a GO-space $L$ if the elements of the pair are the immediate neighbors of each other  with respect to the order of $L$. By $\infty_L$ and $-\infty_L$ we denote the maximum and minimum of $L$ in any ordered compactification. In notations and terminology of general topological nature we will follow \cite{Eng}.
\par\bigskip
\section{Study}

\par\bigskip\noindent
Before we begin our study let us reflect on the structure of  groups of  monotonic autohomeomorphisms on  a GO-space. First, if a GO-space under a discussion is a LOTS, then "monotonic autohomeomorphism" is equivalent to a  "monotonic bijection". In this study, all maps are monotonic autohomeomorphsims.  It is obvious that the composition of two increasing  maps is increasing. If $g$ is decreasing, then $f\circ g$ and $g\circ f$ have the character opposite to that of $f$. We will use the following statement without reference.

\par\bigskip\noindent
\begin{lem}\label{lem:increasingclopen}
Let $L$ be a GO-space. Then the set of decreasing (increasing) autohomeomorphisms on $L$ is clopen in  $M_p(L)$.
\end{lem}
\begin{proof}
Let $f$ be a decreasing autohomeomorphism and $a<b$ Since $f(a)>f(b)$ there exist open intervals $I$ and $J$ containing $f(a)$ and $f(b)$, respectively, such that $I$ is strictly to the right of $J$. Then $U=\{g: g(a)\in I\ and\ g(b)\in J\}$ contains $f$ and does not not contain any increasing autohomeomorphisms. Hence, the set of all decreasing autohomeomorphisms on $L$ is open in $M_p(L)$. Similarly, the set of all increasing autohomeomorphisms is open in $M_p(L)$. Since there are no other elements in $M_p(L)$ the conclusion follows.
\end{proof}

\par\bigskip\noindent
Note that any monotonic autohomeomorphism maps extremities to extremities. This and Lemma \ref{lem:increasingclopen} imply that sets $U(\infty_L; \{y\})$ and $V(-\infty_L;\{y\})$ are open in $M_p(L)$ for any $y\in L$. We will use this fact implicitly throughout the paper.  To ignite our study let us start  with the following positive observation.

\par\bigskip\noindent
\begin{lem}\label{lem:composition}
Let $L$ be a GO-space. Then the operation of function composition is a continuous map from $M_p(L)\times M_p(L)$ to $M_p(L)$.
\end{lem}
\begin{proof}
Fix $f,g\in M_p(L)$. Let $z= f(g(x))$ and $y=g(x)$. Let $W_{f\circ g}$ be an arbitrary neighborhood $f\circ g$.  Our goal is to find open neighborhoods $V_f$ and $U_g$ of $f$ and $g$, respectively, such that $f_1\circ g_1\in W_{f\circ g}$ whenever $f_1\in V_f$ and $g_1\in U_g$. For our argument we will assume that $f$ and $g$ are increasing. Other variations are treated using very similar arguments.
The structure of basic neighborhoods and Lemma \ref{lem:increasingclopen} allow us to assume that $W_{f\circ g}$ is in the form $W_{f\circ g}(x;I)= \{h: h(x) \in I\}$ for some fixed $x\in L$ and a convex open set $I\subset L$. 
We have the following three cases:
\begin{description}
	\item[\rm Case ($z$ is isolated )] Then, $x$ and $y$ are isolated too. Put, $U_f = U_f(x; \{y\})$ and $V_g=V_g(y;\{z\})$.

	\item[\rm Case ($z$ is isolated on one side only)] Without loss of generality, $z$ is a limit point of $\{x\in L: x<z\}$. Hence, $\{x\in L : x > z\}$ is clopen in $L$. Therefore, there exists $z' < z$ such that $[z', z] \subset I$. Since $f$ is onto, there exists $y'\in L$ such that  $f(y')=z'$. Since $f$ is increasing, $y'<y$. Since $f$ is a monotonic homeomorphism, $f([y', y]) = [z',z]$. Put $V_f = V_f(y, y'; I, I)$. Clearly, $V_f$ is an open neighborhood of $f$. 
Put $U_g = U_g(x; (y', y])$. To show that the selected neighborhoods are as desired, pick $f_1\in V_f$ and $g_1\in U_g$. We have $g_1(x)\in (y',y]$. Since $f_1$ is monotonic, we have $f_1(g_1(x))$ is between $f_1(y')$ and $f_1(y)$. By the definition of $V_f$ and convexity of $I$, $f_1(g_1(x))$  is in $I$. Hence,  $f_1(g_1(x))\in W_{f\circ g}$.

	\item[\rm Case ($z$ is a limit point on both sides)]  Since $z$ is limit on both sides, so are $x$ and $y$. Fix $y_1,y_2\in L$ such that $y_1<y<y_2$ and $f(y_1),f(y_2)\in I$ Next, fix $x_1,x_2\in L$ such that $x_1<x<x_2$ and $g(x_1),g(x_2)\in (y_1, y_2)$. By monotonicity, $g(x)$ is between $g(x_1)$ and $g(x_2)$ while $f(y)$ is between $f(y_1)$ and $f(y_2)$. Put $U_g=U_g(x_1,x_2; (y_1,y_2),(y_1,y_2)\}$ and $V_f=V_f(y_1,y_2;I,I)$. Clearly, the sets contain $g$ and $f$, respectively. Fix $g_1\in U_g$ and $f_1\in V_f$. Then $g_1(x)$ is between $g_1(x_1)$ and $g_1(x_2)$, and therefore, $g_1(x)\in (y_1,y_2)$. Since $f_1\in V_f$, $f_1((y_1,y_2))\subset I$. Hence, $f_1(g_1(x))\in I$. Therefore, $f\circ g(x)\in I$.
\end{description}
\end{proof}

\par\bigskip\noindent
In connection with our observation, it must be mentioned that the operation of inversion need not be continuous on $M_p(L)$ when $L$ is a GO-space. 

\par\bigskip\noindent
\begin{ex}\label{ex:inversionfails}
The operation of inversion is not continuous on $M_p(A)$, where $A$ is the Alexandroff Arrow space.
\end{ex}
\begin{proof}
First, recall that the arrow $A$ is the set $(0,1]$ endowed with the topology generated by sets in form $(a,b]$. Let $f$ be the identity map. Put $V_{f^{-1}} = \{h^{-1}: h^{-1}(1/2) \in (1/4,1/2]\}$. Let $U_f$ be any open neighborhood of $f$. We may assume that $U_f$ is in form $U_f(x_1=1/2, x_2, ...,x_n; (a, 1/2], I_2, ..., I_n\}$. Clearly, we can can find an increasing autohomemorphism $g\in U_f$ with the property that $g(1/2) < 1/2$. Then $g^{-1}(1/2) >1/2$, meaning that $g^{-1}\not \in V_{f^{-1}}$.
\end{proof}

\par\bigskip\noindent
\begin{lem}\label{lem:necessary}
Let $L$ be a GO-space. If the operation of inversion is continuous on $M_p(L)$, then any set in form $W(x;\{y\})$ is open in $M_p(L)$ whenever $x$ is isolated from at least one side.
\end{lem}
\begin{proof}
Fix an arbitrary $x\in L$, which is isolated from at least one side and any $y\in L$. Without loss of generality we may assume that $x$ isolated from the right. 
To show that $W=W(x; \{y\})$ is open in $M_p(L)$, fix an arbitrary $f\in W$. We may assume that $f$ is increasing. We need to find an open neighborhood of $f$ which is  a subset of $W$. Since $f\in W$, we have $f(x) = y$. Since $x$ is isolated on the right, the set $I=\{z\in L: z\leq x\}$ is open. Put $V_{f^{-1}} = \{h^{-1}: h^{-1}(y) \in I\}$. Since the operator of inversion is continuous on $M_p(L)$, there exists an open neighborhood $U_f$ of $f$ such that $(U_f)^{-1}\subset  V_{f^{-1}}$. We may assume that there exist $x_2,...,x_n$ and open convex mutually disjoint sets $I_1,...,I_n\subset L$ such that $U_f$ is the set of all increasing functions of  $U(x_1=x, x_2, ... x_n; I_1, .. I_n)$. Since $x$ is isolated from the right and $f$ is increasing we may assume that $\max I_1 = y$. It remains to show that $U_f\subset W$. For this fix $g\in U_f$. We already know that $g(x)\leq y$. Assume that $g(x)<y$. Since $g$ is increasing, we conclude that  $g^{-1}(y)>x$. This  contradicts the fact that $(U_f)^{-1}\subset  V_{f^{-1}}$. Therefore, $g(x) =y$. Hence, $U_f\subset W$. 
\end{proof}

\par\bigskip\noindent
\begin{lem}\label{lem:sufficient}
Let $L$ be a GO-space. If any set in form $W(x;\{y\})$ is open in $M_p(L)$ whenever $x$ is isolated from at least one side, then  the operation of inversion is continuous in $M_p(L)$.
\end{lem}
\begin{proof}
\par\bigskip\noindent
Let $W_{f^{-1}} = \{h^{-1}: h^{-1}(y)\in I\}$ for some fixed $y\in L$ and convex open $I$ in $L$. Let $x=f^{-1}(y)$. We need to find $U_f$ an open neighborhood of $f$ such that $g^{-1}\in W_{f^{-1}}$ for every $g\in U_f$.  We have the following two cases:
\begin{description}
	\item[\rm Case ($x$ is isolated on at least one side )] Then $U_f = \{h: h(x) \in \{y\}\}$ is open by hypothesis. Clearly, $f\in U_f$. Pick any $g\in U_f$. Then, $g^{-1}(y)=x\in I$. Hence, $g^{-1}\in W_{f^{-1}}$.
	\item[\rm Case ($x$ is a limit point on both sides)]  
Fix $x_1,x_2, x_1', x_2'\in I$ such that $x_1 <x_1' < x< x_2'< x_2$. Let $(y_1,y_2)= f((x_1,x_2))$. Without loss of generality, $f$ is increasing. Put $U_f = U_f(x_1',x_2'; (y_1,y), (y,y_2))$. Clearly, $f\in U_f$. Fix $h \in U_f$. 
Then $h(x_1')< y< h(x_2')$. Hence, $h^{-1}\in (x_1',x_2')\subset I$.
\end{description}
\end{proof}

\par\bigskip\noindent
Lemmas \ref{lem:composition}, \ref{lem:necessary}, and \ref{lem:sufficient} form the following criterion.

\par\bigskip\noindent
\begin{thm}\label{thm:criterion}
Let $L$ be a GO-space. Then, $M_p(L)$ is a topological group if and only if  any set in form $W(x;\{y\})$ is open in $M_p(L)$ whenever $x$ is isolated from at least one side.
\end{thm}

\par\bigskip\noindent
We already know that the space of monotonic autohomeomorphisms of the Alexandroff Arrow is not  a topological group. Let us next discuss some positive cases. Firstly, it was proved by Sorin in \cite{Sor} that $M_p(L)$ is a topological group if $L$ is a LOTS. Sorin stated his result for the space of order-preserving bijections but the argument is valid for the space of  all monotonic bijections. To derive Sorin's result from our criterion we need the following lemma.
\par\bigskip\noindent
\begin{lem}\label{lem:gap}
Let $L$ be a GO-space,  $\{a_l, a_r\}\subset L$ a gap,  and $b\in L$. Then $U(f(a_l; \{b\})$ and $V(f(a_r; \{b\})$  are open in $M_p(L)$.
\end{lem}
\begin{proof} Put $S=\{f: f(a_l)=b, \ f\ is\  increasing\}$. By Lemma \ref{lem:increasingclopen}, it suffices to shows that $S$ is open. If there is no increasing  $f$ in $M_p(L)$ that maps $a_l$ to $b$, then $S$ is empty, and therefore, open. Otherwise, fix  $f\in S$. There exists a gap $\{b_l, b_r\}$ such that $b=b_l$, $f(a_l)=b$ and $f(a_r)=b_r$. We have 
$$
U = \{f: f(a_l) = b, f(a_r)=b_r\}=\{f: f(a_l) \in [-\infty_L,b]_L, f(a_r)\in [b_r, \infty_L)_L\}
$$
Since $\{b, b_r\}$ is a gap, the intervals in the rightmost side of the three-sided equality are open in $L$. Hence, $U$ is open in $M_p(L)$. Next, observe that $S=U$.
\end{proof}

\par\bigskip\noindent
Lemma \ref{lem:gap} and Theorem \ref{thm:criterion} imply Sorin's result.
\par\bigskip\noindent
\begin{cor} (Sorin \cite{Sor})
Let $L$ be a LOTS. Then $M_p(L)$ is a topological group.
\end{cor}

\par\bigskip\noindent
Let $x$ be isolated from one side  in a GO-space $L$ and  let $x$ not belong to a gap. Suppose that there exists an open neighborhood  $I$ of $x$ in $L$ that has no other points of this kind other than $x$. Let us show that $U=U(x;\{y\})$ is open in $M_p(L)$ for any $y$. For this pick $f\in U$. Since $f$ is a homeomorphism,  $y$ is also isolated   from one side and is not a member of a gap. Since $f$ is a homeomorphism, $J=f(I)$ is an open neighborhood of $y$ that has no points with the properties of $y$ other than $y$ itself . Put $V_f=\{h:h(x)\in J\}$. Since the inclusion $h(x)\in J$ implies that $h(x)=y$, we conclude that $V_f=U$. This observation leads to the following statement.

\par\bigskip\noindent
\begin{thm}
Let $L$ be a GO-space. If $L$ is a disjoint union of clopen sets each of which is a LOTS, then $M_p(L)$ is a topological group.
\end{thm}

\par\bigskip\noindent
It is obvious that if two LOTS are order-isomorphic, then their spaces of monotonic bijections are homeomorphic and even topologically isomorphic (as topological groups). If, however, two LOTS are simply homeomorphic, their spaces of monotonic bijections need not be homeomorphic. For example, $M_p(\mathbb N)$ contains only the identity map, while $M_p(\mathbb Z)$ is infinite. This is a rather cheap example but it gives a route for exploration. 
Let $\mathcal M_p[L] = \{M_p(L'): L'\ is\ a\ GO-space\ and\ is\ homeomorphic\ to\ L\}$.

\par\bigskip\noindent
\begin{problem}
Identify nice classes  $\mathcal P$ of GO-spaces within which two GO-spaces $L$ and $L'$ are homemorphic if and only if $\mathcal M_p[L]= \mathcal M_p[L']$.
\end{problem}

\par\bigskip\noindent
\par

\end{document}